\input amstex
\magnification=1200
\documentstyle{amsppt}
\NoRunningHeads
\NoBlackBoxes
\define\sltwo{\operatorname{\frak s\frak l}(2,\Bbb C)}
\define\Diff{\operatorname{Diff}}
\define\Vect{\operatorname{Vect}}
\define\HS{\operatorname{\Cal H\Cal S}}
\define\mdhbar{mod\,\hbar}
\title On the infinite-dimensional hidden symmetries. III.
$q_R$--conformal symmetries at $q_R\to\infty$ and Berezin-Karasev-Maslov
asymptotic quantization of $C^\infty(S^1)$
\endtitle
\author Denis V. Juriev
\endauthor
\affil\eightpoint\rm
ul.Miklukho-Maklaya 20-180, Moscow 117437 Russia.\linebreak
E-mail: denis\@juriev.msk.ru\linebreak
\tenpoint\linebreak
\ \linebreak
February, 02, 1997\linebreak
E-version: funct-an/9702002
\ \linebreak
\ \linebreak
\ \linebreak
\endaffil
\abstract
The relations between the infinite dimensional geometry of
$q_R$--con\-for\-mal symmetries at $q_R\to\infty$, Berezin quantization
of the Lobachevski{\v\i} plane and Karasev-Maslov asymptotic quantization are
explicated. Various aspects of the ``approximate'' representation theory are
discussed.
\endabstract
\endtopmatter
\document\ \newline
\ \newline
\ \newline
This short paper being the continuation of the previous parts [1] belongs
to the series of articles supplemental to [2], and also lies in lines of the
general ideology exposed in the review [3]. The main purpose of the activity,
which has its origin and motivation presumably in the author's researches [4]
on the experimentally mathematical aspects of interactively
controlled systems (i.e. the controlled systems, in which the control is
coupled with unknown or uncompletely known uncontrolled feedbacks) and
applications, is to explicate the essentially {\sl infinite-dimensional\/}
aspects of the hidden symmetries, which appear in the representation theory of
the finite dimensional Lie algebras and related algebraic structures. The
present series is organized as a sequence of topics, which illustrate this
basic idea on the simple and tame examples without superfluous difficulties
and details as well as in the series [2] but from a bit more geometric point
of view.
\pagebreak

\head Topic Six: $q_R$--conformal symmetries at $q_R\to\infty$ and
Berezin-Karasev-Maslov asymptotic quantization on $C^\infty(S^1)$
\endhead

The main result of this topic, which may be regarded as an illustration to
the book [5], is the following theorem.

\proclaim{Theorem} Let us realize the Lobachevski{\v\i}-Berezin
$C^*$--algebra [6-8] by the boun\-ded operators in the unitarizable
Verma module $V_h$ over the Lie algebra $\sltwo$ [8]. Then the lifting of
the projective $\HS$--pseudorepresentation of the group $\Diff_+(S^1)$ in
$V_h$ [1:Topic 1] to its action in the  Lobachevski{\v\i}-Berezin algebra
realizes an asymptotic Karasev-Maslov-type quantization [5] (see also [9]) of
$C^\infty(S^1)$ [supplied with the natural action of the group $\Diff_+(S^1)$]
at $h\to\frac12$.
\endproclaim

Let us comment the statement of the theorem.

In the Poincare realization of the Lobachevskii plane (the realization in the
unit disk) the Lobachevskii metric may be written as
$w=q_R^{-1}\,dzd\bar{z}/(1-|z|^2)^2$;
and one can construct the $C^*$--algebra (Lobachevski{\v\i}-Berezin algebra)
[6-8], which may be considered as a quantization of such metric [6], namely,
let us consider two variables $t$ and $t^*$, which obey the following
commutation relations:
$[tt^*,t^*t]=0$, $[t,t^*]=q_R(1-tt^*)(1-t^*t)$ (or in an equivalent form
$[ss^*,s^*s]=0$, $[s,s^*]=(1-q_Rss^*)(1-q_Rs^*s)$, where $s=(q_R)^{-1/2}t$);
one may realize such variables by bounded
operators in the Verma module over $\sltwo$ of the weight
$h=\frac{q_R^{-1}+1}2$ [8];
if such Verma module is realized in polynomials of one complex variable $z$
and the action of $\sltwo$ has the form $L_{-1}=z$, $L_0=z\partial_z+h$,
$L_1=z(\partial_z)^2+2h\partial_z$, then the variables $t$ and $t^*$ are
represented by tensor operators $D=\partial_z$ and $F=z/(z\partial_z+2h)$.
These operators are bounded if
$q_R>0$ and therefore one can construct a Banach algebra
generated by them and obeying the prescribed commutation relations;
the structure of $C^*$--algebra is rather obvious: an involution $*$ is
defined on generators in a natural way, because the corresponding tensor
operators are conjugate to each other.

The ``classical'' case corresponds to the limit $q_R\to 0$. However, it is
interesting to consider another limit transition as $q_R\to\infty$. In this
case the Lobachevski{\v\i}-Berezin algebra is reduced to the algebra
$C^{\infty}(S^1)$ and the variable $t^*$ is identified with $t^{-1}$.
However, the algebra $C^\infty(S^1)$ possesses the group $\Diff_+(S^1)$ of
the orientation preserving diffeomorphisms of a circle as a group of
symmetries. This property is ``weakly'' conserved for finite $q_R$. Namely,
the group $\Diff_+(S^1)$ admits projective $\HS$--pseudorepresentations
in the Verma modules $V_h$ over $\sltwo$ [1]; $\HS$--pseudorepresentation
means the representation up to the Hilbert-Schmidt operators; the prefix
"pseudo" is motivated by the analogy with the pseudodifferential calculus
[10]. The infinitesimal versions of the projective
$\HS$--pseudorepresentations
of the group $\Diff_+(S^1)$, the $\HS$--projective representations of its
Lie algebra $\Vect(S^1)$, the algebra of vector fields on a circle, were
considered in [2:Topic 10]. The generators of $\Vect^{\Bbb C}(S^1)$ in
these $\HS$--projective representations are defined by the $q_R$--conformal
symmetries [8] (the tensor operators of spin 2) and have the following form
$$L_k=(\xi+(k+1)h)\partial_z^k\quad (k\ge0),\quad
L_{-k}=z^k\frac{\xi+(k+1)h}{(\xi+2h)\ldots(\xi+2h+k-1)}\quad (k\ge 1),$$
where $\xi=z\partial_z$.

However, it is rather interesting to investigate the asymptotics of the
$\HS$--pseudorepresentations of $\Diff_+(S^1)$ as $q_R\to\infty$ or
$h\to\frac12$. First of all, let us lift the projective
$\HS$--pseudorepresentation of $\Diff_+(S^1)$ in the Verma module $V_h$ to
the Lobachevski{\v\i}-Berezin algebra as well as the corresponding
$\HS$--projective representation of $\Vect(S^1)$. The analysis of the
explicit formulas for the generators of the Lobachevski{\v\i}-Berezin algebra
and $q_R$--conformal symmetries allows to state that the lifting of
the $\HS$--projective representation of $\Vect^{\Bbb C}(S^1)$ realizes
a representation $\mdhbar$ ($\hbar=2h-1=q_R^{-1}$) of this
Lie algebra in sense of [5,9], moreover, this representation is
one by derivatives of the Lobachevski{\v\i}-Berezin algebra
$\mdhbar$. The theorem is just the global version of this
statement.

Note that $C^\infty(S^1)$ is not supplied by the structure of the Poisson
algebra so our case slightly differs from one considered by M.V.Karasev and
V.P.Maslov, who asymptotically quantize Poisson manifolds. However, the
formal difference is subtle and immaterial.

The theorem may be combined with the results of [1:Topic 2].

\remark{Remark 1} The composed representations of the Witt isotopic pair
[1:Topic 2] in the Verma modules $V_h$ over the Lie algebra $\sltwo$
may be lifted to its representation $\mdhbar$ ($\hbar=2h-1$)
in the Lobachevski{\v\i}-Berezin algebra.
\endremark

It is interesting to ``globalize'' this representation $\mdhbar$ of the
Witt isotopic pair.

The theorem may be interpreted in terms of the nonlinear geometric algebra
[11]. Recall that the projective $\HS$--pseudorepresentations of the group
$G=\Diff_+(S^1)$ in the Verma modules $V_h$ supply $G$ with by the noncanonical
groupuscular structures $\Cal G_h$ [1:Topic 1], which are not odular ones.

\remark{Remark 2} The noncanonical groupuscular structures $\Cal G_h$ on
$G=\Diff_+(S^1)$ define the asymptotically canonical groupuscular structure
$\mdhbar$ ($\hbar=2h-1$).
\endremark

Note that groupuscular structures (as well as general loopuscular structures)
naturally appear in the formalism of the nonlinear Poisson brackets [5] and
differential geometry [11] so their asymptotic considerations are of interest.
Such considerations should combine the nonlinear geometric algebra [11] with
the asymptotic algebraic and differential geometries [12].

\remark\nofrills{Questions:}\newline
(1) What groupuscular structure on $G=\Diff_+(S^1)$ appears from the family
$\Cal G_h$ of noncanonical groupuscular structures $\mdhbar^2$
($\hbar=2h-1$)?\newline
(2) Does the family of projective $\HS$--pseudorepresentations of the
group $\Diff_+(S^1)$, which constitute the asymptotic projective
representation $\mdhbar$ of $\Diff_+(S^1)$, define any hypergroup deformation
[13;5:App.II] of this group?
(3) Does the asymptotic representation of the Lie algebra $\Vect(S^1)$ admit
an extension to the asymptotic representation of the Lie algebra
$\operatorname{DOP}_{[\cdot,\cdot]}(S^1)$ of all differential operators on
the circle?
\endremark

Let us formulate some {\bf conclusions}: (1) the relations between
$q_R$--conformal symmetries at $q_R\to\infty$, the Berezin quantization of
the Lobachevski{\v\i} plane and the Karasev-Maslov asymptotic quantization
were explicated, (2) two approaches to the ``approximate'' representation
theory, namely, one of the $\HS$--pseudorepresentations and one of the
asymptotic representations $\mdhbar$ appeared as being closely and
nontrivially connected in the considered case.

\Refs
\roster
\item" [1]" Juriev D., On the infinite-dimensional hidden symmetries. I,II.
E-prints: funct-an/9612004, funct-an/9701009.
\item" [2]" Juriev D., Topics in hidden symmetries. I-V. E-prints:
hep-th/9405050, q-alg/9610026, q-alg/9611003, q-alg/9611019, funct-an/9611003.
\item" [3]" Juriev D.V., An excursus into the inverse problem of
representation theory [in Russian]. Report RCMPI-95/04 (1995) [e-version:
mp\_arc/96-477].
\item" [4]" Juriev D.V., Octonions and binocular mobilevision [in Russian].
Fundam.Prikl.Matem., to appear; Belavkin-Kolokoltsov watch-dog effects in
interactively controlled stochactic dynamical videosystems [in Russian].
Teor.Matem.Fiz. 106(2) (1996) 333-352 [English transl.: Theor.Math.Phys. 106
(1996) 276-290]; On the description of a class of physical interactive
information systems [in Russian]. Report RCMPI/96-05 (1996) [e-version:
mp\_arc/96-459]; Droems: experimental mathematics, informatics and infinite
dimensional geometry [in Russian]. Report RCMPI/96-05$^+$ (1996).
\item" [5]" Karasev M.V., Maslov V.P., Nonlinear Poisson brackets: geometry
and quantization. Amer.Math.Soc., Providence, RI, 1991.
\item" [6]" Berezin F.A., Quantization in complex symmetric spaces [in
Russian]. Izvestiya AN SSSR. Ser.matem. 39(2) (1975) 363-402.
\item" [7]" Klimek S., Lesniewski A., Quantum Riemann surfaces. I.
Commun.Math.Phys. 146 (1992) 103-122.
\item" [8]" Juriev D., Complex projective geometry and quantum projective
field theory [in Russian]. Teor.Matem.Fiz. 101(3) (1994) 331-348 [English
transl.: Theor.Math.Phys. 101 (1994) 1387-1403];
\item" [9]" Karasev M.V., Maslov V.P., Asymptotic and geometric quantization
[in Russian]. Uspekhi Matem.Nauk 39(6) (1984) 115-173.
\item"[10]" Taylor M., Pseudo-differential operators. B., 1974; Duistermaat J.,
Fourier integral operators. N.Y., 1973; Treves F., Introduction to
pseudo-differential and Fourier integral operators. N.Y., 1980.
\item"[11]" Sabinin L.V., Differential geometry and quasigroups [in Russian].
In ``Current problems of geometry. To the 60th Anniversary of
Acad.Yu.G.Reshetnyak''. Trans.Inst.Math. Siberian Branch Soviet Acad.Sci.,
Novosibirsk, 1984, v.14, pp.208-221; On the nonlinear geometric algebra
[in Russian]. In ``Webs and quasigroups''. Kalinin [Tver'], 1988, pp.32-37;
Mikheev P.O., Sabinin L.V., Quasigroups and differential geometry. In
``Quasigroups and loops. Theory and applications''. Berlin, Heldermann Verlag,
1990, P.357-430.
\item"[12]" Dubnov V.L., Maslov V.P., Naza{\v\i}kinski{\v\i} V.E.,
The complex Lagrangian germ and the canonical operator. Russian J.Math.Phys.
3 (1995) 141-180.
\item"[13]" Litvinov G.L., Hypergroups and hypergroup algebras [in Russian].
Current Probl.Math., Fundam.Directions, V.26, Moscow, VINITI, 1985, P.57-106;
Levitan B.M., Theory of the generalized shift operators [in Russian]. Moscow,
Nauka, 1973; Gurevich D.I., Generalized shift operators on Lie groups [in
Russian]. Izvestiya AN Armyanskoi SSR. 18(4) (1983) 305-317.
\endroster
\endRefs
\enddocument